\documentclass[twoside,11pt]{article}
\usepackage[english]{babel}
\usepackage{amsmath}
\usepackage{amsthm}
\usepackage{cases}
\usepackage{amsfonts}
\usepackage{amssymb}
\usepackage{graphicx}
\usepackage{colortbl,dcolumn}
\usepackage{float}
\usepackage[utf8]{inputenc}
\usepackage[T1]{fontenc}

\usepackage{paralist}  
\allowdisplaybreaks[4]

\newtheorem{theorem}{\bf Theorem}[section]

\newtheorem{corollary}{\bf Corollary}[section]

\newtheorem{remark}{\bf Remark}[section]

\numberwithin{equation}{section}

\begin{document}
	\title{{\Large On self-similar singular solutions to a vorticity stretching equation}
		\footnotetext{\small E-mail addresses: dudp954@nenu.edu.cn (DD),\  lijy645@nenu.edu.cn (JL).} }
	
	\author{{Dapeng Du, Jingyu Li and Xinyue Shi}\\[2mm]
		\small\it School of Mathematics and Statistics, Northeast Normal University,\\
		\small\it   Changchun, 130024, P.~R.~China
	}

	\date{}
	
	\maketitle
	
	\begin{quote}
		\small \textbf{Abstract}:
We consider the following model equation:
\begin{equation}
	\omega_{t} = Z_{11}\omega\,\omega ,
\end{equation}
where \begin{equation}
	Z_{11} = \partial_{11}\Delta^{-1}
\end{equation}
is a Calderon-Zygmond operator.  We get the existence of self-similar singular solutions  with a special form.  The main difficulty is the degeneracy of the operator $Z_{11}$ that is overcome by the spectral uncertainty principle.  We also show that the solution to this model blows up in finite time if the initial datum is compactly supported and has a positive integral.
		

	\end{quote}

\section{Introduction}

The incompressible Euler equations describe the motion of an inviscid incompressible flow.  Whether this system in 3D has singular solutions is one of the most important problems in contemporary PDE theory.  This problem is also the first major step towards the possible singular solutions to the 3D incompressible Navier-Stokes equations,  one of the famous seven millennium  prize problems.

The first result on the singular solutions of incompressible Euler equations is a 1D model constructed by Constantin, Lax and Majda \cite{CLM}.  Thanks to the nice properties of Hilbert transform, they solved the model explicitly.  Since then there are rich developments.  In 2019, an impressive progress was made by Elgindi \cite{El-annals}.  He showed that the solutions to the original incompressible Euler equation blow up in finite time for some initial $ C^{1, \alpha}$ velocity with $\alpha$ small enough.  The idea could trace back to \cite{KS}.  Roughly speaking,  the main thing is that when $\alpha$ goes to zero,  the limiting system of the Euler equations has several nice self-similar singular solutions. The numerical studies by Hou and his collaborators \cite{LH} are also very interesting.
For detailed literature, one may refer to the survey papers \cite{El-survey, Elgindi-survey25}.

The first 2D model appeared in \cite{CS} (with a sign change $w=-w$) has the form
\begin{equation}\label{1.1}
	\omega_{t} = Z_{11}\omega\,\omega.
\end{equation}
For a wide class of negative initial datum,  Constantin and Sun \cite{CS} constructed global regular solutions to this model.  When the whole space is changed to an ellipse,  some explicit singular solutions was constructed in \cite{Du-Euler}.   Very recently, Bianchini and Elgindi \cite{BE} showed that the equation \eqref{1.1}  blows up in finite time  for a large class of initial data.

Here comes our theorem.
\noindent
\begin{theorem}\label{thm-1}
Assume A is a bounded measurable set in $\mathbb{R}^{2}$, then the equation \eqref{1.1} admits a self-similar singular solution int the form
\begin{equation}\label{1.2}
	\omega(x,t)=\frac{1}{T-t} Q \left( x \right) \text{ for } T > 0,
\end{equation}
where
\begin{equation}\label{1.3}
	Q \in L^{2}(\mathbb{R}^{2}), \; Q|_{\mathbb{R}^{2} \backslash A}=0 , \; Z_{11}Q|_{A}=1.
\end{equation}

\end{theorem}

\begin{remark}
This is the first general result on the existence of self-similar singular solutions for the 2D nonlocal models.
\end{remark}

The key in the proof of Theorem \ref{thm-1} is the spectral uncertainty principle.   Roughly speaking,  this principle is the quantitative version of the following property of Fourier transform:  if both the function and its Fourier transform are compactly supported, then this function is identically zero.

It was suggested in \cite{Du-special} that a direct construction of singular solutions may be the way toward singularities for the incompressible Euler and Navier-Stokes equations.  Using the same spectral uncertainty principle, we get a blow-up result.

\begin{theorem}\label{thm-2}
Assume that the equation \eqref{1.1} is equipped with initial data $\omega(x,0) = \omega_{0}$.  If $\omega_{0} \in C^{\alpha}$ with $\alpha>0$, $\omega_{0}$ is compactly supported in $\mathbb{R}^{2}$ with  $\int_{\mathbb{R}^{2}}\omega_{0} dx>0$,  then the solution to \eqref{1.1} blows up in finite time.
\end{theorem}

\begin{remark}  The blow-up criterion of Theorem 1.2 in \cite{BE} is

\begin{itemize}

\item[1)] there exists a $W_{2}\in L^{2}(\mathbb{R}^{2})$ with $W_{2} \ge 0$;

\item[2)]  $Z_{11}W_{2} \, \omega_{0} \ge 0$;

\item[3)]  $-\infty<\int \left(\frac{Z_{11}W_{2}\,\omega_{0}}{W_{2}}\right)W_{2}\,dx $ and  $0 < \int Z_{11}W_{2}\, \omega_{0}\,dx< \infty$.

\end{itemize}

There seems no apparent  way to deduce our theorem from the results in \cite{BE}. The blow-up criterion in our theorem is also simpler.
\end{remark}

\begin{remark}
A series of models  to the incompressible Euler and Navier-Stokes equations were proposed in \cite{Du-Euler}.  It was also conjectured that the singular solutions to 3D incompressible Navier-Stokes equations  generically are largely fluctuated.
\end{remark}

\section{Self-similar singular solution}
\label{sec-2} In this section we will prove  theorem \ref{thm-1}.  First we give a version of spectral uncertainty principle.

\begin{theorem}[Logvinenko-Sereda theorem]\label{thm-3}	
Let $f\in L^{2}(\mathbb{R}^{n}) $ and  $\hat{f}$ be the Fourier transform of f. If  $supp(\hat{f}) \subset A \Subset \mathbb{R}^{n} $,  $B\subset\mathbb{R}^{n}$, $N(B) > 0$,  then  $\|f\|_{L^{2}(B)}\ge C^{-1}\|f\|_{L^{2}(\mathbb{R}^{n})}$ where $C$ is a constant only depending on $A$, $B$ and $n$.
\end{theorem}

\begin{corollary}\label{cor}
	Assume f $\in$ $L^{2}(\mathbb{R}^{n})$,  supp $f \subset A \Subset \mathbb{R}^{n}$,  $B \subset \mathbb{R}^{n}$,  $N(B) \ne 0$.  Then there exists  C=C (A,B,n), such that  $\|\hat{f}\|_{L^{2}(B)} \ge C^{-1}\|\hat{f}\|_{L^{2}(\mathbb{R}^{n})}$,  where $\hat{f}$ is the Fourier transform of f.
\end{corollary}

\begin{proof}	
This proof follows directly from a straightforward application of the Logvinenko-Sereda theorem.
\end{proof}
	
Plugging \eqref{1.2} into \eqref{1.1} implies the following equation for the profile $Q$:
\begin{equation}\label{2.1}
	Z_{11}Q \, Q=Q.
\end{equation}
	
Now we are in a position to prove Theorem \ref{thm-1}.
	
\begin{proof}[Proof of Theorem \ref{thm-1}]
We note that
\begin{equation}
    \eqref{1.2} \iff (Z_{11} \, Q - 1) \, Q = 0.
\end{equation}
Thus for any subset $A \subset \mathbb{R}^2$,  solving \eqref{2.1} reduces to finding $Q$ such that
\begin{equation}
\begin{cases}
1=Z_{11}Q,\quad x \in A ,\\
Q=0, \qquad \; \,    x \notin A  .
\end{cases}
\end{equation}
For $\phi \in L^{2}(a)$ define
\begin{equation}
   \tilde{\phi}=\begin{cases}
	  	\phi ,\quad x \in A ,\\
	  	0    ,\quad x \notin A,
\end{cases}
\end{equation}
and
\begin{equation}  	
 L : L^{2}(A) \longrightarrow L^{2}(A), \   L\phi = Z_{11}\tilde{\phi}|_{A}.
\end{equation}
Next we will prove that $L$ is coercive.
Define
\begin{equation}
	\Omega = \left\{k > \frac{\lambda_{1}}{\lambda_{2}} > \frac{1}{k} |\lambda \ne 0 \right\}.
\end{equation}
Then
\begin{equation}
\begin{split}
	\left\langle L \,\phi,\phi \right\rangle_{L^{2}(A)} =  \left\langle L \,\tilde{\phi},\tilde{\phi} \right\rangle_{L^{2}(R^{2})}
	= \int_{\mathbb{R}^2} F(Z_{11}  \tilde{\phi}) \, F(\tilde{\phi})d\lambda &= \int_{\mathbb{R}^2} \frac{\lambda_1^2}{\lambda_1^2 + \lambda_2^2} |F(\tilde{\phi}) |^2 d\lambda \\
	&\ge \int_{\Omega}\frac{1}{1+k^{2}} | F(\tilde{\phi}) |^2 d\lambda.
\end{split}\end{equation}
From  Corollary \ref{cor} we have
\begin{equation}
	\int_{\Omega}\frac{1}{1+k^{2}} | F(\tilde{\phi}) |^2 d\lambda
	 \ge \frac{1}{1+k^{2}}C_{(\Omega,A))}\int_{R^{2}}| F(\tilde{\phi}) |^2 d\lambda.
\end{equation}
Thus,
\begin{equation}\label{2.9}
	\langle L  \phi,\phi \rangle_{L^{2}(A)} \ge \delta\:\|\phi\|_{L^{2}(A)}^{2} \text{ for }  \delta> 0.
\end{equation}
This implies   $L$ is coercive.

We proceed to show that $L$ is one-to-one and its inverse is bounded.

1) One-to-one. Assume $\phi \in L^{2}(A)$ and $L\,\phi=0$. Define
\[
	\tilde{\phi}=\begin{cases}
		\phi, \quad x\in  A,\\
		0, \quad x\notin A.
	\end{cases}
\]
Then we have
\[
	\langle L \phi , \phi \rangle_{L^{2}(A)} = \langle Z_{11}\,\tilde{\phi}, \tilde{\phi}\rangle_{L^{2}(R^{2})}=0,
\]
which indicates
\[
	\tilde{\phi}\equiv0.
\]
Hence $\phi\equiv0$ and  $L$ is one-to-one.

2) The inverse is bounded. Using \eqref{2.9} we get
\[
\langle L  \phi , \phi \rangle_{L^{2}(A)}=\langle Z_{11}  \tilde{\phi} , \tilde{\phi} \rangle_{L^{2}(R^{2})}\\
\ge \delta \|\phi\|_{L^{2}(A)},
\]
which gives
\[
\|\phi\|_{L^{2}(A)}\le \delta^{-1}  \|\phi\|_{L^{2}(A)} \, \|L \phi \|_{L^{2}(A)}.
\]
It then follows that
\[
\|\phi\|_{L^{2}(A)} \le \delta^{-1}  \|L \phi\|_{L^{2}(A)}.
\]
Hence the inverse of $L$ is bounded.

Next we show $L$ is onto. Assume the contrary. Then $\exists \, \phi_{0} \notin L(L^{2}(A))$.  Since the inverse of $L$ is bounded,  the space $L(L^{2}(A))$ is closed. Denote $\phi_{1}$ the projection of $\phi_{0}$ on $L(L^{2}(A))$. Then we have
\[
\langle \phi_{0}-\phi_{1}, L(L^{2}(A)) \rangle = 0.
\]
Thus,
\[
\langle L(\phi_{0}-\phi_{1}),\phi \rangle = 0, \forall \, \phi \in L^{2}(A),
\]
which means $L(\phi_{0}-\phi_{1})=0$.  Since $L$ is one-to-one, we get $\phi_{0}=\phi_{1}$, which is a contradiction. Therefore $L$ is onto.

After proving $L$ is onto,  the proof is essential finished. Let  $\tilde{Q}=L^{-1}\,1$ and
\[
Q=\begin{cases}
	\tilde{Q} \quad x\in A,\\
	0 \quad x\notin A.
\end{cases}
\]
Then \eqref{1.2} is a solution to \eqref{1.1}. The proof of Theorem \ref{thm-1} is complete.
\end{proof}

\section{Blow-up}\label{sec-3}

In this section we prove Theorem \ref{thm-2}.

\begin{proof}[Proof of Theorem \ref{thm-2}]	
First  we have
\begin{equation}
	\int_{\mathbb{R}^{2}} \omega_{t} d\lambda = \int_{\mathbb{R}^{2}}Z_{11} \omega\,\omega d\lambda.
\end{equation}
Let $\hat{\omega}$ be the Fourier transform of $\omega$. Then
\begin{equation}\begin{split}
	\int_{\mathbb{R}^{2}}Z_{11}\omega\,\omega d\lambda
	 =\int_{\mathbb{R}^{2}}\widehat{Z_{11}\omega} \hat{\omega}d\lambda =\int_{\mathbb{R}^{2}}\dfrac{\lambda_{1}^{2}}{\lambda_{1}^{2}+\lambda_{2}^{2}}\hat{\omega}\,\hat{\omega}d\lambda
	 \ge \int_{\Omega}\frac{| \hat{\omega}|^{2}}{1+k^{2}}d\lambda.	
\end{split}\end{equation}
It follows from  Corollary  \ref{cor} that
\begin{equation}
	\frac{\int_{\Omega} |\hat{\omega}|^{2}d\lambda}{\int_{\mathbb{R}^{2}}|\hat{\omega}|^{2}d\lambda} \ge C_{(\Omega,\: A)} > 0.	
\end{equation}
Thus,
\begin{equation}
	\int_{\Omega}\frac{|\hat{\omega}|^{2}}{1+k^{2}}d\lambda \ge \frac{C_{(\Omega,\: A)}}{1+k^{2}}\int_{\mathbb{R}^{2}} \left|\hat{\omega}\right|d\lambda,
\end{equation}
from which we obtain that
\begin{equation}
	\omega_{t} \ge \delta\:\omega^{2} \text{ for some }  \delta  > 0.
\end{equation}
Therefore the solution of equation \eqref{1.1} blows up in finite time under the conditions of Theorem \ref{thm-2}.
\end{proof}

\section*{Acknowledgements}
Dapeng Du would like to thank Professor Hongjie Dong and Professor Zhifei Zhang for valuable discussions.


\begin{thebibliography}{99}
	
\bibitem{BE} Roberta Bianchini and Tarek M. Elgindi, Finite-time singularity formation for scalar stretching equations, Nonlinearity 38 (2025), no. 7, Paper No. 075003, 11 pp.


\bibitem{CLM} P. Constantin, P. D. Lax, and A. Majda, A simple one-dimensional model for the three-dimensional vorticity equation, Comm. Pure Appl. Math. 38 (1985), no. 6, 715-724.


\bibitem{CS} P. Constantin and W. Sun, Remarks on Oldroyd-B and related complex fluid models, Commun. Math. Sci. 10 (2012), no. 1, 33-73.

\bibitem{El-survey} Theodore D. Drivas and T. M. Elgindi, Singularity formation in the incompressible Euler equation in finite and infinite time, EMS Surv. Math. Sci. 10 (2023), no. 1, 1-100.

\bibitem{Du-special} Dapeng Du, Special solutions for nonlinear partial differential equations from physics (in Chinese), Journal of Northeast Normal University, 52 (2020), no. 1, 1-3.


\bibitem{Du-Euler} Dapeng Du, On some model equations of Euler and Navier-Stokes equations, Chinese Ann. Math. Ser. B 42 (2021), no. 2, 281-290.


\bibitem{El-annals} T. M. Elgindi,  Finite-time singularity formation for $C^{1,\alpha}$ solutions to the incompressible Euler equations on $\mathbb{R}^{3}$,  Ann. Math. (2) 194 (2021), no. 3, 647-727.


\bibitem{Elgindi-survey25}  T. M. Elgindi, Dynamics of ideal fluid flows, arXiv: 2511.16254v1.


\bibitem{HQWW} De Huang, Xiang Qin, Xiuyuan Wang, and Dongyi Wei, Self-Similar finite-time blowups with smooth profiles of the generalized Constantin-Lax-Majda model,  Arch. Ration. Mech. Anal. 248 (2024), no. 2, Paper No. 22, 65 pp.


\bibitem{KS} A. Kiselev and V. Sverak,  Small scale creation for solutions of the incompressible two-dimensional Euler equation, Ann. Math. (2) 180 (2014), no. 3, 1205-1220.


\bibitem{LH} G. Luo and T. Y. Hou, Potential singular solutions of the 3D incompressible Euler equations, arXiv: 1310.0497v2.

\end{thebibliography}
\end{document}